\documentclass{article}

\usepackage{amsmath}
\usepackage{amssymb}
\usepackage{amsthm}

\usepackage{times}

\newtheorem{lemma}{Lemma}
\newtheorem{theorem}[lemma]{Theorem}
\newtheorem{definition}[lemma]{Definition}
\newtheorem{claim}{Claim}

 \newenvironment{claimproof}{\begin{proof}}{\end{proof}}

\author{Henry Towsner}
\title{A Model Theoretic Proof of Szemer\'edi's Theorem}
\date{\today}
\begin{document}
\maketitle

\begin{abstract}
  We present a short proof of Szemer\'edi's Theorem using a dynamical system enriched by ideas from model theory.  The resulting proof contains features reminiscent of proofs based on both ergodic theory and on hypergraph regularity.
\end{abstract}

\section{Introduction}

Szemer\'edi's Theorem states:
\begin{theorem}
  For any $\delta>0$ and any $k$, there is an $n$ such that whenever $N\geq n$, $A\subseteq [1,N]$, and $|A|\geq\delta N$, there exists an $a$ and a $d$ such that
\[a,a+d,a+2d,\ldots,a+(k-1)d\in A.\]
\end{theorem}

Szemer\'edi's original proof \cite{szemeredi:MR0369312} used graph theoretic methods, in particular the Szemer\'edi Regularity Lemma \cite{szemeredi:MR540024}.  Shortly after, Furstenberg gave a different proof \cite{furstenberg77,furstenberg:MR670131}, based on a correspondence argument which translates the problem into one in ergodic theory.  Beginning with a new proof given by Gowers \cite{gowers:MR1844079}, a number of new proofs have been developed in the last decade.  (Tao counts a total of roughly sixteen different proofs \cite{taoblogmanyproofs}.)

Hrushovski has recently used a stronger correspondence-type argument \cite{hrushovski} to make progress on a similar combinatorial problem (the so-called non-commutative Freiman conjecture).  In this paper, we use Hrushovski's method to give a short proof of Szemer\'edi's theorem.

The proof here bears a similarity to proofs based on hypergraph regularity, such as \cite{gowers:MR2373376,nagle:MR2198495,rodl:MR2069663,tao:MR2259060}; in particular the proof is very similar to the infinitary regularity-like arguments introduced by Tao \cite{tao07} and used by Austin to prove both Szemer\'edi's Theorem \cite{austinmultiszem} and generalizations \cite{austinhalesjewitt}.  Indeed, this proof was inspired by noticing that the use of ``wide types'' (countable intersections of definable sets of positive measure) in Hrushovski's arguments was analogous to the use of the regularity lemma in finitary arguments.  (In fact, Hrushovski essentially sketches a proof of the $k=3$ case of Szemer\'edi's Theorem in \cite{hrushovski}; however his arguments depend on stability theoretic methods which don't seem to generalize to higher $k$.  This seems related to the fact that stability implies $3$-amalgamation, but not $4$-amalgamation.)

The methods here are also reminiscent of those used by Tao to prove the ``diagonal ergodic theorem'' \cite{tao08Norm}, and especially to our infinitary reformulation of that proof \cite{towsner:MR2529651}.  This paper might shed light on the connection between that method and the technique of ``pleasant extensions'' used by Austin \cite{austinnonconventional}.

We are grateful to the members of UCLA's reading seminar on \cite{hrushovski}: Matthias Aschenbrenner, Isaac Goldbring, Greg Hjorth, Terence Tao, and Anush Tserunyan.

\section{A Correspondence Principal}
Suppose that Szemer\'edi's Theorem fails; that is, for some $\delta>0$ and every $n$, there is an $A_n\subseteq [1,n]$ with $|A_n|\geq \delta n$ such that $A_n$ contains no $k$-term arithmetic progression.  We must first make a technical adjustment: we view the sets $A_n$ as subsets of the group $[1,2n+1]$.  Then $|A_n|\geq\delta n/2-\epsilon$ where $\epsilon\rightarrow 0$ as $n\rightarrow\infty$, and $A_n$ contains no $k$-term arithmetic progressions in this group.\footnote{Without this modification, $A_n$ might contain no arithmetic progressions as a set of integers, but contain arithmetic progressions as a subset of the \emph{group}, since the group contains progressions which ``wrap around'': for instance, $3,10,4$ is an arithmetic progression in $[1,13]$, but does not correspond to an arithmetic progression in the integers.  Expanding the group by adding a ``dead zone'' disjoint from $A_n$ is a standard way of avoiding this problem.}

We extend the language of groups with a predicate symbol $A$ and the following additional class of formulas:
\begin{itemize}
\item Whenever $\alpha_1,\ldots,\alpha_k,\gamma$ is a sequence of rationals, $\vec x$ a tuple, and $B_1(\vec x,\vec y_1),\ldots,B_k(\vec x,\vec y_k)$ a sequence of formulas, $\int \sum_{i\leq k}\alpha_i\cdot B_idm(\vec x)>\gamma$ is a formula with free variables $\vec y_1,\ldots,\vec y_k$
\end{itemize}
We let $([1,2n+1],A_n)$ be models, interpreting the symbol $A$ by $A_n$ and 
\[\left[\int \sum_{i\leq k}\alpha_i\cdot\phi_idm(\vec x)>\gamma\right](\vec p_1,\ldots,\vec p_k)\Leftrightarrow \sum_{i\leq k}\alpha_i\frac{|\{\vec x\mid \phi_i(\vec x,\vec p_i)\}|}{n^{|\vec x|}}>\gamma\]
We write $m_{\vec x}(\phi)$ as an abbreviation for $\int 1\cdot\phi dm(\vec x)$, and sometimes omit $\vec x$ when it is clear from context.

Form an ultraproduct of those groups $[1,2n+1]$ such that $2n+1$ is prime.  We obtain a model $(G,+,A)$.  By transfer, the formula
\[\exists a,d (a\in A\wedge a+d\in A\wedge\cdots a+d+\cdots+d\in A)\]
is false.

Observe that for any countable set $M$ and any $n$, the sets of $n$-tuples definable with parameters from $M$ form an algebra of internal sets, and using the Loeb measure construction, we may extend the internal counting measure on this model to a measure $\mu^n$ on the $\sigma$-algebra of Borel sets generated from these definable sets (for basic facts about this construction, see \cite{Goldblatt}).  The measures $\mu^n$ satisfy Fubini's Theorem \cite{Keisler1,Keisler2}; that is,
\[\int f d\mu^n=\iint f d\mu^{n_0}d\mu^{n_1}\]
where $n_0+n_1=n$.

\begin{lemma}
  If $(G,+,A)\vDash m(\phi)>\gamma$ then $\mu^{|\vec x|}(\{\vec x\mid \phi(\vec x)\})\geq\gamma$.

If $(G,+,A)\vDash \neg m(\phi)>\gamma$ then $\mu^{|\vec x|}(\{\vec x\mid \phi(\vec x)\})\leq\gamma$.
\end{lemma}
\begin{proof}
For the first part, since $(G,+,A)\vDash m(\phi)>\gamma$, for almost every $n$, $([1,2n+1],+,A_{2n+1})\vDash m(\phi)>\gamma$, and therefore
\[\frac{|\{\vec x\mid\phi(\vec x)\}|}{(2n+1)^{|\vec x|}}>\gamma\]
holds in $([1,2n+1],A_{2n+1})$.  (Here $\phi$ may contain parameters.)

But since this holds for almost every $n$, by transfer
\[\frac{|\{\vec x\mid\phi(\vec x)\}|}{(2n+1)^{|\vec x|}}>\gamma\]
holds in $(G,+,A)$ where the inequality is between nonstandard rational numbers.  This means that
\[\mu^{|\vec x|}(\{\vec x\mid\phi(\vec x)\})=st(\frac{|\{\vec x\mid\phi(\vec x)\}|}{(2n+1)^{|\vec x|}})\geq\gamma.\]

For the second part, suppose $(G,+,A)\vDash \neg m(\phi)>\gamma$.  Then for almost every $n$, $([1,2n+1],+,A_{2n+1})\vDash \neg m(\phi)>\gamma$, and therefore
\[\frac{|\{\vec x\mid\phi(\vec x)\}|}{(2n+1)^{|\vec x|}}\leq\gamma\]
holds in $([1,2n+1],A_{2n+1})$. But since this holds for almost every $n$, by transfer,
\[\frac{|\{\vec x\mid\phi(\vec x)\}|}{(2n+1)^{|\vec x|}}\leq\gamma\]
holds in $(G,+,A)$, so also
\[\mu^{|\vec x|}(\{\vec x\mid\phi(\vec x)\})=st(\frac{|\{\vec x\mid\phi(\vec x)\}|}{(2n+1)^{|\vec x|}})\leq\gamma.\]
\end{proof}
Since both $m$ and $\mu$ are additive, this extends immediately to the corresponding integrals.

As a notational convenience, let us write $\vec x$ for a sequence $x_1,\ldots,x_k$, and $\vec x_{\overline i}$ for the sequence $x_1,\ldots,x_{i-1},x_{i+1},\ldots,x_k$ and $\vec x_{\overline {i,k}}$ for the sequence $x_1,\ldots,x_{i-1},x_{i+1},\ldots,x_{k-1}$.

Observe that for any $b$, the function $a\mapsto a+b$ is a definable bijection, as is $a\mapsto k\cdot a$ for any integer $k$.  For $r\leq k$, we define formulas $A_r$ on $x_{\overline r}$.  When $r<k$, we define
\[A_r(x_{\overline r}):\Leftrightarrow \sum_{i<k,i\neq r}i\cdot x_i+r\cdot\left(x_k-\sum_{i<k,i\neq r} x_i\right)\in A.\]
We define
\[A_k(x_{\overline k}):\Leftrightarrow\sum_{i<k}i\cdot x_i\in A.\]
Note that for any $x_{\overline {k-1,k}}$, $\mu(\{x_{k-1}\mid x_{\overline k}\in A_k\})=\mu(A)$, so by the Fubini property of these measures, $\mu^{k-1}(A_k)=\mu(A)$.

Define $\sigma(\vec x_{\overline k}):=\sum_{i<k}x_i$; it is easy to see that for any $\vec a_{\overline {i,k}}$, the function $a_i\mapsto\sigma(\vec a_{\overline k})$ is a bijection.  Whenever $A_k(\vec a_{\overline k})$, also $A_i(\vec a_{\overline {i,k}},\sigma(\vec a_{\overline k}))$ for each $i<k$.  In particular, if we define $\hat A:=\{\vec a_{\overline k}\mid A_k(\vec a_{\overline k})\wedge\forall i<k A_i(\vec a_{\overline {i,k}},\sigma(\vec a_{\overline k}))\}$, we have $\hat A=A_k$, and so $\mu^{k-1}(\hat A)=\mu(A)>0$.

In the next section, we will show that, under these conditions, $\mu^k(\bigcap_{i\leq k}A_i)>0$.  First, however, we show that this is enough to prove Szemer\'edi's Theorem.  If $\mu^k(\bigcap_{i\leq k}A_i)>0$, we may find $\vec a\in\bigcap_{i\leq k}A_i$ such that $a_k\neq\sum_{i<k}a_i$.  Setting $a:=\sum_{i<k}i\cdot a_i$ and $d:=a_k-\sum_{i<k}a_i$, we have $a+id\in A$ for $i\in[0,k-1]$.  Therefore $a,a+d,\ldots,a+(k-1)d$ is a $k$-term arithmetic progression in $A$.

 This contradicts the construction of the model, which in turn means that the initial assumption that the sets $A_N\subseteq[1,N]$ exist must fail.  Therefore for every $\delta$, there is an $N$ such that for every $A_N\subseteq[1,N]$ with $|A_N|\geq\delta n$, $A_N$ contains an arithmetic progression of length $k$.

\section{Amalgamation}

\begin{definition}
Let $M\subseteq G$ be a set, let $n$ be a positive integer, and let $I\subseteq [1,n]$ be given.  We define $\mathcal{B}_{n,I}(M)$ to be the Boolean algebra of subsets of $G^n$ of the form
\[\{(x_1,\ldots,x_n)\mid\phi(\{x_i\}_{i\in I})\}\]
where $\phi$ is a formula with parameters from $M$.

We write $S\choose k$ for the collection of subsets of $S$ with cardinality $k$.  When $k\leq n$, we define $\mathcal{B}_{n,k}(M)=\bigvee_{I\in{[1,n]\choose k}}\mathcal{B}_{n,I}(M)$.

We also define $\mathcal{B}_{n,<I}(M)=\bigvee_{J\in{I\choose |I|-1}}(M)$.

If $\mathcal{B}$ is any Boolean algebra, we write $\mathcal{B}^\sigma$ for the $\sigma$-algebra generated by $\mathcal{B}$.
\end{definition}

We equate formulas with the sets they define, so we will also speak of $\mathcal{B}$ as being a Boolean algebra of formulas. 

\begin{lemma}
Let $\mathcal{B},\mathcal{B}_0,\mathcal{B}_1$ be Boolean algebras with $\mathcal{B}\subseteq\mathcal{B}_0\cap\mathcal{B}_1$.  Suppose there is an elementary submodel $M$ such that:
\begin{itemize}
  \item Every parameter in every formula in $\mathcal{B}_1$ belongs to $M$,
  \item If $\phi(\vec x,\vec a)\in\mathcal{B}_0$, $\vec b\in M$, and $|\vec a|=|\vec b|$, then $\phi(\vec x,\vec b)\in\mathcal{B}$.
\end{itemize} 
Then for any $f\in L^2(\mathcal{B}_1^\sigma)$, $||E(f\mid\mathcal{B}^\sigma)-E(f\mid\mathcal{B}_0^\sigma)||=0$.
\end{lemma}
An illustrating case is when $\mathcal{B}=\mathcal{B}_{n,k}(M)$, $\mathcal{B}_0=\mathcal{B}_{n,k}(M\cup N)$, and $\mathcal{B}_1=\mathcal{B}_{n,k+1}(M)$. 
\begin{proof} 
Suppose not.  Then setting $\epsilon:=||f-E(f\mid\mathcal{B}^\sigma)||_{L^2}$ and $\delta:=||f-E(f\mid\mathcal{D}_0^\sigma)||_{L^2}$, we must have $\delta<\epsilon$.  For some $\beta_1,\ldots,\beta_m$ and $A_1(\vec x,\vec b_1),\ldots,A_m(\vec x,\vec b_m)\in\mathcal{B}_1$ with each $\vec b_i\in M$, we have $||f-\sum_{i\leq m}\beta_i\chi_{A_i}(\vec x,\vec b_i)||_{L^2}<(\epsilon-\delta)/4$.

Since $||f-E(f\mid\mathcal{B}_0^\sigma)||=\delta$, there are $\alpha_1,\ldots,\alpha_n$ and $D_1(\vec x,\vec a_1),\ldots,D_n(\vec x,\vec\alpha_n)\in\mathcal{D}_0$ with $||f-\sum_{i\leq n}\alpha_i\chi_{D_i}(\vec x,b_i)||_{L^2}<\epsilon-(\epsilon-\delta)/2$, and therefore
\[||\sum_{i\leq m}\beta_i\chi_{A_i}(\vec x,\vec b_i)-\sum_{i\leq n}\alpha_i\chi_{D_i}(\vec x,b_i)||_{L^2}<\epsilon-3(\epsilon-\delta)/4.\]

This means the formula
\[\exists\vec y_1,\dots,\vec y_m\neg(||\sum_{i\leq m}\beta_i\chi_{A_i}(\vec x,\vec b_i)-\sum_{i\leq n}\alpha_i\chi_{D_i}(\vec x,\vec y_i)||^2_{L^2}>\left[\epsilon-3(\epsilon-\delta)/4\right]^2)\]
is satisfied (where, to view this as a formula, we expand the norm into an integral of sums of definable formulas).

This is a formula with parameters from $M$, so by the elementarity of $M$, there are witnesses $\vec a'_1,\ldots,\vec a'_n$ in $M$ satisfying:
\[||\sum_{i\leq m}\beta_i\chi_{A_i}(\vec x,\vec b_i)-\sum_{i\leq n}\alpha_i\chi_{D_i}(\vec x,\vec a'_i)||_{L^2}\leq\epsilon-3(\epsilon-\delta)/4.\]

It follows that $||f-\sum_{i\leq n}\alpha_i\chi_{D_i}(\vec x,\vec a'_i)||_{L^2}<\epsilon$, and since $\sum_{i\leq n}\alpha_i\chi_{D_i}(\vec x,\vec a'_i)$ is measurable with respect to $\mathcal{D}^\sigma$, this contradicts the assumption that $||f-E(f\mid\mathcal{D}^\sigma)||_{L^2}=\epsilon$.
\end{proof}

\begin{theorem}
Let $n\leq k$ and suppose that for each each $I\subseteq[1,n]$ with $|I|=k$, we have a set $A_I\in\mathcal{B}^\sigma_{n,I}(M)$, and suppose there is a $\delta>0$ such that whenever $B_I\in\mathcal{B}_{n,I}(M)$ and $\mu^n(A_I\setminus B_I)<\delta$ for all $I\in{[1,n]\choose k}$, $\bigcap B_I$ is non-empty.  Then $\mu^n(\bigcap_{I\in{[1,n]\choose k}} A_I)>0$.
\end{theorem}
\begin{proof}
We proceed by main induction on $k$.  When $k=1$, the claim is trivial: we must have $\mu(A_I)>0$ for all $I$, since otherwise we could take $B_I=\emptyset$; then $\mu(\bigcap A_I)=\prod\mu(A_I)>0$.  So we assume that $k>1$ and that whenever $B_I\in\mathcal{B}_{n,I}(M)$ and $\mu^n(A_I\setminus B_I)<\delta$ for all $I$, $\bigcap_{I\in{[1,n]\choose k}} B_I$ is non-empty.  Throughout this proof, the variable $I$ ranges over ${[1,n]\choose k}$.

\begin{claim}
For any $I_0$,
\[\int (\chi_{A_{I_0}}-E(\chi_{A_{I_0}}\mid\mathcal{B}_{n,<I_0}^\sigma(M)))\prod_{I\neq I_0}\chi_{A_I}d\mu^{n}=0.\]
\end{claim}
\begin{claimproof}
When $k=n$, this is trivial since $\prod_{I\neq I_0}\chi_{A_I}$ is an empty product, and therefore equal to $1$.

If $k<n$, we have
\begin{align*}
 &\int (\chi_{A_{I_0}}-E(\chi_{A_{I_0}}\mid\mathcal{B}_{n,<I_0}^\sigma(M)))\prod_{I\neq I_0}\chi_{A_I}d\mu^{n}\\
 =&\int (\chi_{A_{I_0}}-E(\chi_{A_{I_0}}\mid\mathcal{B}_{n,<I_0}^\sigma(M)))\prod_{I\neq I_0}\chi_{A_I}d\mu^{k}(\{x_i\}_{i\in I_0})d\mu^{n-k}(\{x_i\}_{i\not\in I_0}).\\
\end{align*}
Observe that for any choice of $\{a_i\}_{i\not\in I_0}$, $\mathcal{B}_{n,<I_0}(M),\mathcal{B}_{n,<I_0}(M\cup\{a_i\}_{i\not\in I_0}),\mathcal{B}_{n,I_0}(M)$ satisfy the preceding lemma, so
\[||E(\chi_{A_{I_0}}\mid\mathcal{B}_{n,<I_0}^\sigma(M))-E(\chi_{A_{I_0}}\mid\mathcal{B}_{n,<I_0}^\sigma(M\cup\{a_i\}_{i\not\in I_0}))||_{L^2}=0.\]
The function $\prod_{I\neq I_0}\chi_{A_i}(\{x_i\}_{i\in I_0},\{a_i\}_{i\not\in I_0})$ is measurable with respect to $\mathcal{B}_{n,<I_0}^\sigma(M\cup\{\vec a_i\})$.  Combining these two facts, we have
\begin{align*}
& \int (\chi_{A_{I_0}}-E(\chi_{A_{I_0}}\mid\mathcal{B}_{n,<I_0}^\sigma(M)))\prod_{I\neq I_0}\chi_{A_I}d\mu^{k}(\{x_i\}_{i\in I_0})\\
=&\int (\chi_{A_{I_0}}-E(\chi_{A_{I_0}}\mid\mathcal{B}_{n,<I_0}^\sigma(M\cup\{a_i\}_{i\not\in I_0})))\prod_{I\neq I_0}\chi_{A_i}d\mu^k(\{x_i\})\\
=&0
\end{align*}
Since this holds for any $\{a_i\}_{i\not\in I_0}$, the claim follows by integrating over all choices of $\{a_i\}$.
\end{claimproof}

\begin{claim}
For any $I_0$, there is an $A'_{I_0}\in\mathcal{B}_{n,<I_0}^\sigma(M)$ such that:
\begin{itemize}
  \item Whenever $B_I\in\mathcal{B}_{n,I}(M)$ for each $I$, $\mu^n(A_I\setminus B_I)<\delta$ for each $I\neq I_0$, and $\mu^n(A'_{I_0}\setminus B_{I_0})<\delta$, $\bigcap_{I\in{[1,n]\choose k}}B_I$ is non-empty, and
  \item If $\mu^n(A'_{I_0}\cap \bigcap_{I\neq I_0}A_I)>0$, $\mu^n(\bigcap A_I)>0$.
\end{itemize}
\end{claim}
\begin{claimproof}
 Define $A'_{I_0}:=\{\vec x_i\mid E(\chi_{A_{I_0}}\mid\mathcal{B}_{n,<I_0}^\sigma(M))(\vec x_i)>0\}$.  If $\mu^n(A'_{I_0}\cap\bigcap_{I\neq I_0}A_I)>0$ then we have
 \[\int E(\chi_{A_{I_0}}\mid\mathcal{B}_{n,<I_0}^\sigma(M))\prod_{I\neq I_0}\chi_{A_I}d\mu^{n}>0\]
 and by the previous claim, this implies that $\mu(\bigcap A_I)>0$.
 
 Suppose that for each $I$, $B_I\in\mathcal{B}_{n,I}(M)$ with $\mu^n(A_I\setminus B_I)<\delta$ for $I\neq I_0$ and $\mu^n(A'_{I_0}\setminus B_I)<\delta$.  Since
\[\mu^n(A_{I_0}\setminus A'_{I_0})=\int\chi_{A_{I_0}}(1-\chi_{A'_{I_0}})d\mu^n=\int E(\chi_{A_{I_0}}\mid\mathcal{B}^\sigma_{n,<I_0})(1-\chi_{A'_{I_0}})d\mu^n=0,\]
we have $\mu^n(A_{I_0}\setminus B_{I_0})<\delta$ as well, and therefore $\bigcap B_I$ is non-empty.
\end{claimproof}

By applying the previous claim to each $I\in {[1,n]\choose k}$, we may assume for the rest of the proof that for each $I$, $A_I\in\mathcal{B}_{n,<I}^\sigma(M)$.

Fix some finite algebra $\mathcal{B}\subseteq\mathcal{B}_{n,k-1}(M)$ so that for every $I$, $||\chi_{A_I}-E(\chi_{A_I}\mid\mathcal{B})||_{L^2}<\frac{\sqrt{\delta}}{\sqrt{2}({n\choose k}+1)}$ (such a $\mathcal{B}$ exists because there are finitely many $I$ and each $A_I$ is $\mathcal{B}^\sigma_{n,k-1}(M)$-measurable).  For each $I$, set $A^*_I=\{\vec a_i\mid E(\chi_{A_I}\mid\mathcal{B})(\vec a)>\frac{{n\choose k}}{{n\choose k}+1}\}$.

\begin{claim}
For each $I$, $\mu(A_I\setminus A^*_I)\leq \delta/2$
\end{claim}
\begin{claimproof}
$A_I\setminus A^*_I$ is the set of points such that $\chi_{A_I}-E(\chi_{A_I}\mid\mathcal{B})(\vec a)>\frac{1}{{n\choose k}+1}$.  By Chebyshev's inequality, the measure of this set is at most
\[({n\choose k}+1)^2\int (\chi_{A_I}-E(\chi_{A_I}\mid\mathcal{B}))^2d\mu=({n\choose k}+1)^2||\chi_{A_I}-E(\chi_{A_I}\mid\mathcal{B})||_{L^2}^2<\frac{\delta}{2}.\]
\end{claimproof}

\begin{claim}
$\mu(\bigcap_I A_I)\geq\mu(\bigcap_I A^*_I)/\left({n\choose k}+1\right)$. 
\end{claim}
\begin{claimproof}
 For each $I_0$, 
 \begin{align*}
 \mu((A^*_{I_0}\setminus A_{I_0})\cap\bigcap_{I\neq I_0} A^*_I)
 &=\int \chi_{A^*_{I_0}}(1-\chi_{A_{I_0}})\prod_{I\neq I_0}\chi_{A^*_I}d\mu^n\\
 &=\int\chi_{A^*_{I_0}}(1-E(\chi_{A_{I_0}}\mid\mathcal{B}))\prod_{I\neq I_0}\chi_{A^*_I}d\mu^n\\
 &\leq \frac{1}{{n\choose k}+1}\int\prod_{I}\chi_{A^*_I}d\mu^n\\
 &=\frac{1}{{n\choose k}+1}\mu(\bigcap_I A^*_I)
\end{align*}

But then
\[\mu(\bigcap_I A^*_I\setminus \bigcap_I A_I)
\leq \sum_{I_0}\mu((A^*_{I_0}\setminus A_{I_0})\cap\bigcap_{I\neq I_0} A^*_I)
\leq \frac{{n\choose k}}{{n\choose k}+1}\mu(\bigcap_I A^*_I).\]
\end{claimproof}

Each $A^*_I$ may be written in the form $\bigcup_{i\leq r_I}A^*_{I,i}$ where $A^*_{I,i}=\bigcap_{J\in {I\choose k-1}}A^*_{I,i,J}$ and $A^*_{I,i,J}$ is an element of $\mathcal{B}_{n,J}(M)$.  We may assume that if $i\neq i'$ then $A^*_{I,i}\cap A^*_{I,i'}=\emptyset$.

We have
\[\mu(\bigcap_I A^*_I)=\mu(\bigcup_{\vec i\in\prod_I [1,r_I]}\bigcap_I\bigcap_{J\in{I\choose k-1}}A^*_{i_I,J,I}).\]
For each $\vec i\in\prod_I [1,r_I]$, let $D_{\vec i}=\bigcap_I\bigcap_{J\in{I\choose k-1}}A^*_{i_I,J,I}$.  Each $A^*_{I,i_I,J}$ is an element of $\mathcal{B}_{n,J}(M)$, so we may group the components and write $D_{\vec i}=\bigcap_{J\in{n\choose k-1}}D_{\vec i,J}$ where $D_{\vec i,J}=\bigcap_{I\supset J}A^*_{I,i_I,J}$.

Suppose, for a contradiction, that $\mu(\bigcap_I A^*_I)=0$.  Then for every $\vec i\in\prod_I [1,r_I]$, $\mu(D_{\vec i})=\mu(\bigcap_J D_{\vec i,J})=0$.  By the inductive hypothesis, for each $\gamma>0$, there is a collection $B_{\vec i,J}\in\mathcal{B}_{n,J}(M)$ such that $\mu(D_{\vec i,J}\setminus B_{\vec i,J})<\gamma$ and $\bigcap_J B_{\vec i,J}=\emptyset$.  In particular, this holds with $\gamma=\frac{\delta}{2{k\choose k-1}(\prod_I r_I)(\max_I r_I)}$.

For each $I,i\leq r_I, J\subset I$, define
 \[B^*_{I,i,J}=A^*_{I,i,J}\cap\bigcap_{\vec i, i_I=i}\left[B_{\vec i,J}\cup\bigcup_{I'\supseteq J, I'\neq I}\overline{A^*_{I',i_{I'},J}}\right].\]
 
 \begin{claim}
 $\mu(A^*_{I,i,J}\setminus B^*_{I,i,J})\leq \frac{\delta}{2{k\choose k-1}(\max_I r_I)}.$ 
\end{claim}
\begin{claimproof} 
Observe that if $x\in A^*_{I,i,J}\setminus B^*_{I,i,J}$ then for some $\vec i$ with $i_I=i$, $x\not\in B_{\vec i,J}\cup\bigcup_{I'\supseteq J, I'\neq I}\overline{A^*_{I',i_{I'},J}}$.  This means $x\not\in B_{\vec i,J}$ and $x\in\bigcap_{I'\supseteq J}A^*_{I',i_{I'},J}=D_{\vec i,J}$.  So
\[\mu(A^*_{I,i,J}\setminus B^*_{I,i,J})\leq \sum_{\vec i\in\prod_I [1,r_I]}\mu(D_{\vec i,J}\setminus B_{\vec i,J})\leq \frac{\delta}{2{k\choose k-1}(\max_I r_I)}.\]
\end{claimproof}

 Define $B^*_{I}=\bigcup_{i\leq r_I}\bigcap_J B^*_{I,i,J}$.
 
 \begin{claim}
 $\mu(A_I\setminus B^*_{I})<\delta$. 
\end{claim}
\begin{claimproof}
Since $\mu(A_I\setminus A^*_I)<\delta/2$, it suffices to show that $\mu(A^*_I\setminus B^*_{I})<\delta/2$.

 \begin{align*} 
 \mu(A^*_I\setminus \bigcup_i\bigcap_J B^*_{I,i,J})
&=\mu\left(\bigcup_{i}\bigcap_{J}A^*_{I,i,J}\setminus\bigcup_i\bigcap_J B^*_{I,i,J}\right)\\
&\leq\sum_{i\leq r_I}\mu\left(\bigcap_J A^*_{I,i,J}\setminus\bigcap_J B^*_{I,i,J}\right)\\
&\leq\sum_{i\leq r_I}\sum_{J}\mu(A^*_{I,i,J}\setminus B^*_{I,i,J})\\
&\leq r_I\cdot {k\choose k-1}\cdot \frac{\delta}{2{k\choose k-1}(\max_I r_I)}\\
&\leq\delta/2
\end{align*}
\end{claimproof}

\begin{claim} 
\[\bigcap_{I}B^*_{I}\subseteq\bigcup_{\vec i}\bigcap_J B_{\vec i,J}.\]
\end{claim}
\begin{claimproof} 
Suppose $x\in\bigcap_{I}B^*_{I}=\bigcap_I\bigcup_{i\leq r_I}\bigcap_J B^*_{I,i,J}$.  Then for each $I$, there is an $i_I\leq r_I$ such that $x\in\bigcap_J B^*_{I,i_I,J}$.  Since $B^*_{I,i_I,J}\subseteq A^*_{I,i_I,J}$, for each $I$ and $J\subset I$, $x\in A^*_{I,i_I,J}$.

For any $J$, let $I\supset J$.  Then
\[x\in B^*_{I,i_I,J}=A^*_{I,i_I,J}\cap\bigcap_{\vec i',i'_I=i_I}\left[B_{\vec i,J}\cup\bigcup_{I'\supseteq J, I'\neq I}\overline{A^*_{I',i_{I'},J}}\right].\]
In particular, $x\in \left[B_{\vec i,J}\cup\bigcup_{I'\supseteq J, I'\neq I}(\overline{A^*_{I',i_{I'},J}})\right]$ for the particular $\vec i$ we have chosen.  Since $x\in A^*_{I,i_{I'},J}$ for each $I'\supset J$, it must be that $x\in B_{\vec i,J}$.  This holds for any $J$, so $x\in \bigcap_J B_{\vec i,J}$.
\end{claimproof}

From our assumption, $\bigcap_{I}B^*_{I}$ is non-empty, and therefore there is some $\vec i$ such that $\bigcap_J B_{\vec i,J}$.  But this leads to a contradiction, so it must be that $\mu(\bigcap_I A^*_I)>0$, and therefore, as we have shown, $\mu(\bigcap_I A_I)\geq\frac{1}{{n\choose k}+1}\mu(\bigcap_I A^*_I)>0$.
\end{proof}

\bibliographystyle{plain}
\bibliography{../../Bibliographies/main}

\end{document}